\documentclass{proc-l}

\usepackage{amssymb}
\usepackage{amsthm}

\newtheorem{theorem}{Theorem}[section]
\newtheorem{lemma}[theorem]{Lemma}

\theoremstyle{definition}
\newtheorem{definition}[theorem]{Definition}

\newtheorem{proposition}[theorem]{Proposition}
\newtheorem{claim}[theorem]{Claim}

\theoremstyle{remark}

\newcommand{\force}{\mbox{$\Vdash$}}

\newcommand{\ga}{\alpha}
\newcommand{\gb}{\beta}
\newcommand{\grg}{\gamma}
\newcommand{\gd}{\delta}
\newcommand{\gre}{\varepsilon}

\newcommand{\gk}{\kappa}
\newcommand{\gl}{\lambda}
\newcommand{\gm}{\mu}

\newcommand{\gt}{\tau}

\newcommand{\go}{\omega}

\newcommand{\ha}{\aleph}
\newcommand{\rng}{\hbox{rng}}

\newcommand{\mG}{\hbox{{\bf G}}}
\newcommand{\mH}{\hbox{{\bf H}}}

\newcommand{\mP}{\hbox{{\bf P}}}

\newcommand{\mS}{\hbox{{\bf S}}}

\newcommand{\mPt}{\mbox{$\bf \tilde P$}}

\newcommand{\Dt}{\mbox{$\tilde D$}}

\newcommand {\va}{\hbox{\bf a}}

\newcommand {\vi}{\hbox{\bf i}}

\newcommand{\jh}{\mbox{$\hat j$}}

\newcommand{\cf}{\hbox{cf}}

\newcommand{\los}{\mbox{$L[0^{\sharp}]$}}

\newcommand{\rarrow}{{\rightarrow}}

\newcommand{\restrict}{{\upharpoonright}}

\numberwithin{equation}{section}

\begin {document}
\title {$0^{\sharp}$ and Elementary End Extensions of $V_{\gk}$}
\author {Amir Leshem}
\address{Institute of Mathematics, Hebrew University, Jerusalem
, Israel}
\curraddr{Circuit and Systems, Faculty of Information Technology and 
Systems, Mekelweg 4, 2628CD Delft, The Netherlands}
\email{leshem@cas.et.tudelft.nl}
\date {}
\commby{Carl G. Jockusch}
\begin{abstract}
In this paper we prove that if $\gk$ is a cardinal in
$\los$, then there is an inner model $M$ such that
$M \models (V_{\gk},\in)$ has no elementary end
extension.   In particular if $0^{\sharp}$ exists
then weak compactness is never downwards absolute.
We complement the result with a lemma stating that any cardinal greater than 
$\ha_1$ of uncountable cofinality in $\los$ is Mahlo in every strict inner 
model of $\los$.
\end{abstract}
\maketitle
\bibliographystyle{plain}
\section{Introduction}
In this paper we consider the question of existence of elementary
end extensions of models of the form $(V_{\gk},\in)$.
\begin{definition}
\begin {itemize}
\item[1.]
Let $({\mathbb E}_M,\prec_M )$ denote the structure of
all non-trivial elementary end
extensions of M, with $A\prec_M B$ iff B is an elementary
end extension of A.
\item[2.]
Let $({\mathbb E}_M^{\hbox{wf}},\prec_M )$ denote the structure of
all non-trivial {\em well founded} elementary end
extensions of M, with $A\prec_M B$ iff B is an elementary
end extension of A.
\end {itemize}
\end{definition}
Several results regarding the existence of elements in
${\mathbb E}_M$ were proved by Keisler, Silver and Morley.
\begin{theorem}[Keisler, Morley] Let M be a model of ZFC,
$cof(On^M)=\go$. Then ${\mathbb E}_M \neq \emptyset$.
\end{theorem}
\begin{theorem}[Keisler,Silver] Let $M=(V_{\gk},\in)$ be a model of
ZFC, where $\gk$ is weakly compact cardinal. Then for every $S \subseteq M$
${\mathbb E}_{(V_{\gk},\in,S)}^{\hbox{wf}} \neq \emptyset$.
\end{theorem}
Villaveces \cite{Vi96a}, \cite{Vi96b}
has proved several other results regarding the existence
of elementary end extensions of $V_{\gk}$.
\begin{theorem}[Villaveces]
The theory ``ZFC + GCH + $\exists \lambda (\lambda$ measurable) +
$\forall \kappa[\kappa$ inaccessible not weakly compact
$\rarrow  \exists$ transitive
$M_\kappa \models
ZFC$ such that $o(M)=\kappa$ and ${\mathbb E}^{\hbox{wf}}_M=\emptyset$]" is
consistent relative to the theory ``ZFC + $\exists \lambda (\lambda$
measurable) + the weakly compact cardinals are cofinal in On''.
\end{theorem}
He also proved that the property
${\mathbb E}_{V_{\gk}}^{\hbox{wf}} \neq \emptyset$, is not preserved in certain
generic extensions by destroying a weakly compact cardinal.
In this paper we consider the  problem of downwards absoluteness of
the existence of well founded elementary end extensions of $V_{\gk}$. We 
prove the following :
\begin{theorem}
\label{main_theorem}
If $0^{\sharp}$ exists then for every cardinal $\gk$ there
is an inner model M such that
\begin{equation}
M \models {\mathbb E}_{V_{\gk}}=\emptyset.
\end{equation}
\end{theorem}
In particular weak compactness is never downwards absolute,
once we have $0^{\sharp}$ in the universe. On the other hand we will
prove that any cardinal with uncountable cofinality is Mahlo in any
strict inner model of $\los$. I would like to thank the referee for pointing 
out an inaccuracy in the formulation of lemma \ref{Mahlo_lemma} and for 
asking the question at the end of the paper.
\section{Main Theorem}
In this section we prove theorem \ref{main_theorem}.
Let $\gk$ be a cardinal.
Since we assume that $0^{\sharp}$ exists we can
construct our model inside the inner model $L[0^{\sharp}]$ . Note
that since $\gk$ is a cardinal in $V$ it remains a cardinal
in  $L[0^{\sharp}]$, and hence it is weakly compact in $L$.
Our model will be a generic extension of $L$, such that we will be
able to construct a generic object inside $L[0^{\sharp}]$.
The basic idea will be to construct a generic Suslin tree
and then to code it. For the construction of the Suslin
tree we  will follow Kunen's construction \cite{Ku78}, while the
coding will use Levy collapse of certain $L$ cardinals.
 Then we will obtain the generic filter inside $L[0^{\sharp}]$.
%A word about notation, since $L[G]$ can be used either as a
%generic extension of $L$ or as the constructible universe relative to $G$,
%we emphasize that in the discussion below except $L[0^{\sharp}]$
%all the extensions are {\em generic} extensions.

The following theorem by Kunen gives us the forcing for
generating the Suslin tree.
\begin{theorem}
\label{Kun}
Let $\gk$ be a weakly compact cardinal and $P_{\gk}$ be the forcing for 
adding a Cohen subset to
$\gk$. Then $P_{\gk} \simeq R_{\gk} * T_{\gk}$, where
$R_{\gk}$ is a forcing that adds a Suslin tree $T_{\gk}$
to $\gk$, and $T_{\gk}$ is the forcing defined by the tree.
\end{theorem}
Let $\mP$ be the reverse Easton iteration for adding a
Cohen subset to each inaccessible, defined by :
\begin{definition}
\begin{equation}
\mP=(P_{\ga},Q_{\ga} | \ga \in On),
\end{equation}
where
\begin {itemize}
\item[]$P_0=\emptyset$.
\item[] If $\ga$ is not inaccessible then
        $P_{\ga} \force Q_{\ga}=\emptyset$
\item[] If $\ga$ is inaccessible then $Q_{\ga}$ is a
$P_{\ga}$ name for a partial order adding a Cohen subset to $\ga$
i.e.
$P_{\ga} \force Q_{\ga}=\left(2^{<\ga},\subseteq\right)$.
\end {itemize}
Direct limits are taken at inaccessible limits of inaccessibles and
inverse limits otherwise.
\end{definition}
Solovay (see M. Stanley \cite{St88}) proved that the reverse Easton
support iteration for adding Cohen subsets to every
$L$ inaccessible has a generic filter in $L[0^{\sharp}]$, and therefore
our iteration up to $\gk$ has a generic filter as well.

Let $\mG=\left<G_{\ga} | \ga \le \gk \right>$ be $\mP$ generic.
By Kunen's theorem we can interpret $G_{\gk}$ as a pair
$G_{\gk}=\left<T_{\gk},b_{\gk}\right>$ where $T_{\gk}$ is a
$\gk$ Suslin tree and $b_{\gk}$ is a
branch through $T_{\gk}$.

Next we define the forcing used to code the tree $T_{\gk}$.
Let $\mS$ be the Easton supported product of collapsing
of $\ga^{+3}$ to $\ga^{+2}$ defined inside $L$.
\begin{equation}
\mS=\prod \left\{ S_{\ga} :  \ga \hbox{\ is inaccesible \ }\right\}
\end{equation}
where
$S_{\ga}=Coll(\ga^{+2},\ga^{+3})$.

\begin{proposition}
\label{generic_prp}
There is a
 $\mP \times \mS$ generic over $L$, inside $L[0^{\sharp}]$.
\end{proposition}
\begin{proof}
The method of proof of this lemma is
almost identical to the proof of M. Stanley
of Solovay's theorem that there exists a $\mP$ generic filter over $L$
inside $L[0^{\sharp}]$.
We shall build the generic filter by induction on the
Silver indiscernibles. The main point will be taking care that at limits
the generic filter will be the direct limit of the
previously built generic filters.

Let $\left<i_{\ga} : \ga<\gk \right>$ be an increasing enumeration of
the indiscernibles below $\gk$. For any indiscernible $\gl$ the forcing
can be factored as
\begin{equation}
\mP \times \mS =
\left(\mP^{\gl+1}*\mP_{\gl+1}\right) \times \left(\mS^{\gl} \times
\mS_{\gl}\right)
\end{equation}
where $\mP^{\gl}$ is the iteration up to $\gl$, and $\mP_{\gl}$
is the iteration from $\gl$ upwards.
For each $\ga$ we shall define
$\left(G^{i_\ga},H^{i_\ga}\right)$, and then define
$\left(G^{i_\ga+1},H^{i_\ga+1}\right)$ such that
$G^{i_\ga+1} \times H^{i_\ga+1}$ is
$\left(\mP^{i_\ga}*Q_{i_\ga}\right) \times \mS^{i_\ga+1}$ generic
over $L[G^{i_\ga} \times H^{i_\ga}]$.\\
\underline {$i_{0}$ or $i_{\ga+1}$. \ }\\
We have that in $L$ for every indiscernible $\gl$
both $\mP'_{\gl+1}$ and $\mS_{\gl}$ are $\gl^{++}$
closed,  where
\begin{equation}
\mP'_{\gl+1}=\left\{ \gt : \gt
\hbox{\ is a name and \ } \force_{\mP^{\gl+1}}
\gt \in \mPt_{\gl+1} \right\}
\end{equation}
 is the term forcing for
$\mP_{\gl+1}$.
Hence $\mP_{\gl+1} \times \mS_{\gl}$ is $\gl^+$-distributive
over $L^{\mP^{\gl+1} \times \mS^{\gl}}$,
since $\mP^{\gl+1}\times \mS^{\gl}$ is
obviously $\gl^+$-c.c.

By the same argument
$\mP_{i_{\ga}+1}^{i_{\ga+1}}
\times \mS_{i_{\ga}}^{i_{\ga+1}}$ is also $i_{\ga}^+$
distributive.
Let
\begin{equation}
M=L^{\mP^{i_{\ga}+1} \times \mS^{i_{\ga}+1}}.
\end{equation}
Note that each $L$ name for dense subset of
$\mP_{i_{\ga}+1}^{i_{\ga+1}} \times \mS_{i_{\ga}}^{i_{\ga+1}}$
in $M$, belongs to the Skolem hull of the ordinals up to $i_{\ga}$
and finitely many indiscernibles above $i_{\ga+1}$, say
$\{i_{\ga+1},\ldots,i_{\ga+n}\}$. Hence in $\los$ we can represent the
dense subsets of $\mP_{i_{\ga}+1}^{i_{\ga+1}} \times
\mS_{i_{\ga}+1}^{i_{\ga+1}}$
in $M$, by a countable union of families of dense subsets each of size
$i_{\ga}$. Now using the $i_{\ga}^+$ distributivity we can meet each of
these dense subsets. To ensure downwards compatibility we also demand that 
$\left(G^{i_\ga+1},H^{i_\ga+1}\right)$ extends 
$\left(G^{i_\ga},H^{i_\ga}\right)$. 
Finally use the same distributivity argument to define
a generic filter $G(i_{\ga+1})$ for $Q_{i_{\ga+1}}$ over
$L^{\left(\mP^{i_{\ga+1}} \times \mS^{i_{\ga+1}}\right)}$. Again in order to 
ensure extension we demand that  $G(i_{\ga+1})$ extends $G(i_{\ga})$, by
putting a condition forcing it into the generic. Since $\mS$ is not active at 
these stages and using the fact that $\mP$ is a reverse Easton iteration 
this is possible.

\underline{$i_{\ga}$ for $\ga$ limit.}\\
We have built generic objects 
$\left<G^i \times H^i : i<\ga\right>$
for the product up to $\ga$. Now we would like to build a generic filter for
$\mP^{i_\ga} \times \mS^{i^\ga}$. Note that since $i_{\ga}$ is Mahlo
in $L$ we take direct limit. Moreover $\mP^{i_\ga} \times \mS^{i_\ga}$ is
$i_\ga-$c.c.
Define $G^{i_\ga},H^{i_\ga}$ by
\begin{equation}
p\in G^{i_\ga} \hbox{\ iff \ }
\forall \grg<i_\ga p\restrict \grg \in G^{i_\grg}.
\end{equation}
\begin{equation}
s\in H^{i_\ga} \hbox{\ iff \ }
\forall \grg<i_\ga s\restrict \grg \in H^{i_\grg}.
\end{equation}
We prove that $G^{i_\ga} \times H^{i_\ga}$ is
$\mP^{i_\ga} \times \mS^{i_\ga}$ generic over $L$. Suppose that
\\ $D \subseteq \mP^{i_\ga} \times \mS^{i_\ga}$ is dense open.
$D$ belongs to the Skolem hull of finitely many ordinals below
$i_\ga$ $\va=\left<\grg_1,\ldots,\grg_n\right>$ and
finitely many indiscernibles above $\ga$ say
$\vi_n=\left<i_{\ga+1},\ldots,i_{\ga+n}\right>$.
Let $\sup(\va)<i_\gb<i_\ga$.
Define an elementary embedding $j:L \rarrow L$ by
\begin{equation}
j(i_{\grg})= \left\{
\begin{array}{cc}
i_{\grg} & \hbox{\ if \ } \grg<\gb \\
i_{\ga+\gd} & \hbox{\ if \ } \grg=\gb+\gd,0 \le \gd
\end{array} \right.
\end{equation}
Obviously $D \in \rng j$, and $j^{-1}(D)$ is dense open in
$\mP^{\gb} \times \mS^\gb$. Let
$(p',q') \in j^{-1}(D) \cap \left(\mP^{\gb} \times \mS^\gb\right)$.
Since both $p',q'$ are trivial on an end segment we obtain
that
\begin{equation}
j((p',q'))=(p,q)^{\wedge}\left< \emptyset^{Q_\grg \times S_{\grg}} :
\gb \le \grg < \ga\right>.
\end{equation}
Hence by our choice of $\left(G^{i_\ga},H^{i_\ga}\right)$
we obtain that
$j((p',q')) \in \left(G^{i_\ga},H^{i_\ga}\right)$.

Finally we prove that we can find a generic object $G(i_\ga)$
for $Q_{i_\ga}$ over
$L^{\left(G^{i_\ga}\times H^{i_\ga}\right)}$. Define
\begin{equation}
G(i_\ga)=\cup_{\gb<\ga} G(i_\gb).
\end{equation}
Let $D$ be a dense subset of $Q_{i_\ga}$ in
$L^{G^{i_\ga} \times H^{i_\ga}}$.
Let $\Dt$ be a name for $D$ in
$\mP^{i_\ga} \times \mS^{i_\ga}$
Again $\Dt$ is in the Skolem hull of some $i_\gb<i_\ga$
and finitely many indiscernibles
$\vi_n=(i_{\ga+1},\ldots,i_{\ga+n})$.
Define $j:L \rarrow L$ as above. As we have proved
if $(p,q) \in G^{i_\gb} \times H^{i_\gb}$ then
$j(p,q) \in G^{i_\ga} \times H^{i_\ga}$.
Hence the embedding $j$ has a canonical extension to an
embedding
$\jh:L[G^{i_\gb} \times H^{i_\gb}]
\rarrow L[G^{i_\ga}\times H^{i_\ga}]$
defined by
\begin{equation}
\jh(\gt(G^{i_\gb} \times H^{i_\gb}))=
j(\gt)(G^{i_\ga}\times H^{i_\ga}).
\end{equation}
Since $\Dt$ is in $\rng j$ we have $D \in \rng \jh$.
The proof ends as follows: \\
Let
\begin{equation}
p' \in G(i_\gb) \cap \jh^{-1}(D).
\end{equation}
$p'$ exists since by induction hypothesis. 
$G(i_{\gb})$ is $Q_{i_\gb}$ generic, and $\jh^{-1}(D)$
is dense in $Q_{i_\gb}$ by elementarity, and hence 
$\jh(p') \in D$. Since
$p' \in L_{i_\gb}[G^{i_\gb} \times H^{i_\gb}]$
we have
$j(p')=p'$. So
\begin{equation}
p' \in G(i_\gb) \cap D \subseteq G(i_\ga) \cap D.
\end{equation}
\end{proof}

Let $\mG \times \mH$ be $\mP \times \mS$ generic over $L$.
Suppose that $\mH=\left< h_{\ga} |\ga<\gk\right>$ is the
$\mS$ generic filter.
Let $<\cdot,\cdot>$ be a definable pairing function in L, such that for every
$\gb,\grg$, $<\gb,\grg>$ is an L inaccessible.
Since the pairing is definable and $\gk$ is an indiscernible it is closed 
under the pairing function.

Let $T$ be the tree part of $G(\gk)$.
Our final model will be
$N=L[T,\left<h_{\ga} | \ga \in C_T\right>]$
where
\[
C_T=\{\ga|\exists \gb,\grg \left(\ga=<\gb,\grg>  \wedge \gb <_T\grg
\right).
\]
To finish the proof of the theorem we have to prove:
\begin{proposition}
\label{end_prp}
\begin{equation}
N \models ``V_{\gk} \hbox{\ has no elementary end extension".}
\end{equation}
\end{proposition}
\begin{proof}
%First note that $V_{\gk}^N=L_{\gk}[T,\left<h_{\ga} |\ga \in C_T\right>]$.
The proof will be done by a sequence of claims.
\begin{claim}
$N \models T \hbox {\ is Suslin}$.
\end{claim}
\begin{proof}
The claims follows from the fact that the forcing $\mS$ is $\gk$-Knaster in
$L[T]$.
Hence $\mS \times T$ is $\gk$-c.c. in $L[T]$, so especially $T$ is
$\gk$-c.c. in $N'=L[T,\left<h_{\ga}|\ga<\gk\right>]$.
But $N \subseteq N'$ and $\gk^N=\gk^{N'}$, thus $N$ contains no large
anti-chains of $T$ as well.
\end{proof}
\begin{claim}
For every inaccessible $\ga$
\begin{equation}
N \models \ga^{+++L}<\ga^{+++} \iff \ga \in C_T.
\end{equation}
\end{claim}
\begin{proof}
Since for every $\ga \in C_T$ the claim obviously holds, it will be
enough to prove that other cardinals are not collapsed inside
$L[\mG,\left<h_{\ga}|\ga \in C_T \right>]$.
For each $\gm \not \in C_T$ we can even work inside
$L[\mG,\left<h_{\ga}|\ga \neq \gm \right>]$.
However since both forcing
notions $\mP$ and
\[
\mS^{-\gm}=\prod \left\{ S_{\ga} : \ga \neq \gm
\hbox{\ and $\ga$ is inaccesible \ }\right\}
\]
factors nicely, it is obvious that
the only $L$-cardinals collapsed are the triple successors of cardinals
in $C_T$.
\end{proof}

Notice that by the inaccessibility of $\gk$ all the collapsing functions
are inside $V_{\gk}^N$.

Now we finish the proof of proposition \ref{end_prp}.
In $(V_{\gk}^N,\in)$ the tree $T$
is definable by the first order formula:
\[
\gb<_T\grg \iff \exists \ga \left( \ga \hbox{\ is inaccessible\ } \wedge
\ga=<\gb,\grg> \wedge \ga^{+++L}<\ga^{+++}\right).
\]
$(V_{\gk}^N,\in) \models T \hbox {\ is a $\gk$ tree}$,i.e., for every ordinal
$\ga$ $\{x \in T | \hbox{hight}_T(x)=\ga\}$ is a set, and for every ordinal
$\ga$ there is an element of $T$ of hight $\ga$.
Assume that $(M,E)$ is an end extension of  $(V_{\gk}^N,\in)$. 
Let $a$ be a new ordinal in $M$.  In $M$
there is a tree $T'$ which end extends the tree $T$, since $T$ was definable.
By elementarity 
\[
M \models \hbox{there is a branch $b$ in $T'$ of length $a$}. 
\]
Now it follows that
\[
N \models \{x \in b | rk(x)<\gk \} \hbox{ is a branch through \ } T.
\]
Hence any 
end extension of $(V_{\gk}^N,\in)$ will provide a branch through $T$ in $N$.
This is a contradiction since $N \models T \hbox {\ is Suslin}$.
\end{proof}

\section{Mahloness in inner models}
In view of the previous result it is natural to ask whether we can get an 
inner model $M \subsetneq \los$  such that for every inaccessible cardinal 
$\ga \in M$, $(V_{\ga},\in)$ has no well founded elementary end extension.
This turns out to be impossible by the following lemma:
\begin{lemma}
\label{Mahlo_lemma}
Let $\gk>\ha_1$ be a cardinal in $L[0^{\sharp}]$, $\cf(\gk)>\ha_0$, then $\gk$ is
weakly Mahlo in any strictly inner model $M \subsetneq L[0^{\sharp}]$.
Moreover if $\gk$ is a limit cardinal then $\gk$ is strongly Mahlo in every 
$M \subsetneq \los$.
\end{lemma}
\begin{proof}
 The basic idea is to use the covering theorem
to prove that certain cardinals are not
collapsed, in any strict inner model of $\los$. Then we use the covering theorem
again to  prove that actually there must be a stationary set of inaccessibles
below $\gk$.
Let $M \subsetneq \los$ be an inner model.
Let $I=\left\{i_{\ga} | \ga \in On \right\}$ be an increasing
enumeration of Silver's indiscernibles. Then for every $\ga$ such that
$\go<\cf (\ga)$ we have $M \models i_{\ga}^{+L}$ is a cardinal. The proof of 
this uses an idea of Beller \cite{BJW}.
Assume $M\models  i_{\ga}^{+L}$ is not a cardinal. Then 
$|i_{\ga}^{+L}|^M=|i_{\ga}|^M$. By the covering theorem also
$M\models \cf(i_{\ga}^{+L})=\cf(i_{\ga})=|i_{\ga}|$. So in $M$ there is an $f$
$f:i_{\ga} \rarrow i_{\ga}^{+L}$ which maps in an order preserving way a
cofinal subset of $i_{\ga}$ into a cofinal 
subset of $i_{\ga}^{+L}$. Since 
$\los\models \cf(i_{\ga}^{+L})=\go$ choose a cofinal sequence 
(in $i_{\ga}^{+L}$) $\{\gb_n:n<\go\}$ inside $\los$. Now let
$\grg_n$ be the least $\grg$ such that $f(\grg)>\gb_n$.
We obtain that  $\{\grg_n:n<\go\}$ is cofinal in $i_{\ga}$ so  
$\los\models \cf(i_{\ga})=\go$. This contradicts the fact that $i_{\ga}$ has
uncountable cofinality.
Hence every limit of indiscernibles of uncountable
cofinality is a limit cardinal. By the covering theorem it must be a
regular cardinal, so it is weakly inaccessible. Especially any uncountable 
cardinal is weakly inaccessible.

Suppose now that $i_{\ga}$ is not Mahlo in $M$ and $i_{\ga}$ is a limit of 
indiscernibles of uncountable cofinality. Then there is a club 
$C \subseteq i_{\ga}$ consisting of singular cardinals in $M$. By the covering 
theorem (between $L$ and $M$) each element of $C$ is singular in $L$. Hence
$C \cap I = \emptyset$. Hence $\los \models \cf(i_{\ga})=\go$ (since it has
two disjoint clubs through $i_{\ga}$). Therefore if $\cf(i_{\ga})>\go$
and $i_{\ga}$ is a limit of indiscernibles of uncountable cofinality 
it must be Mahlo in any strict inner model. 

If  $\gk$ is also a limit cardinal in $\los$ it is strong limit by GCH. 
Hence it is strong limit in any inner model, so it is strongly Mahlo in $M$. 
\end{proof}

Therefore if  $\gk$ is limit in $\los$ and $\cf(\gk)>\go$, then 
in every inner model there is an inaccessible $\ga<\gk$ such that
${\mathbb E}^{\hbox{wf}}_{(V_{\ga},\gre)} \not = \emptyset$.

A natural question is whether one can have no weakly compacts in a strictly 
inner model of $\los$. We comment that if 
there is a $\gk$ such that $\los \models \gk \rarrow(\go)^{< \go}$ then by a 
result of  Silver \cite{silver66}, any inner model $M$, 
$M \models \gk \rarrow(\go)^{< \go}$, hence there are many ineffable 
cardinals in $M$. Similarly if there is a subtle 
cardinal $\gk$, in $\los$, then obviously $\gk$ is subtle in every inner model
(the definition is $\Pi_1$). Hence there are many large cardinals below it in 
any inner model (e.g., totally indescribables).

However the following question remains open:\\
{\bf Question:} $(ZFC+V=\los)$. Let $M$ be an inner model. Is it consistent 
that $M$ has no weakly compact cardinals ? Is it consistent that for no $\gk$
$M \models \gk \rarrow (\go)^{<\go}$ ?

%\bibliography{logic}

\end{document}